\documentclass[12 pt]{article}

\usepackage{latexsym}
\usepackage{amssymb, amsmath}
\usepackage{amsthm}
\usepackage{enumerate}
\usepackage{geometry}
\usepackage{hyperref}
\usepackage{changebar}
\usepackage{tikz}
\usetikzlibrary{calc}
\usetikzlibrary{decorations.markings, arrows.meta}
\usepackage{todonotes}
\usepackage[normalem]{ulem}

\numberwithin{equation}{section}

\newtheorem{theorem}{Theorem}[section]
\newtheorem{lemma}[theorem]{Lemma}

\newtheorem{proposition}[theorem]{Proposition}

\newtheorem{remark}[theorem]{Remark}

\newcommand{\cA}{\mathcal{A}}

\newcommand{\cK}{\mathcal{K}}

\newcommand{\cC}{\mathcal{C}}

\newcommand{\cM}{\mathcal{M}}
\newcommand{\cB}{\mathcal{B}}

\newcommand{\cS}{\mathcal{S}}

\newcommand{\cQ}{\mathcal{Q}}
\newcommand{\cR}{\mathcal{R}}

\def\un{\underline}

\newcommand\ignore[1]{}

\newcommand{\vf}{\varphi}

\newcommand{\R}{\mathbb{R}}
\newcommand{\ve}{\varepsilon}

\newcommand{\musrb}{\mu_{\mbox{\tiny SRB}}}

\newcommand\md[1]{{\color{orange}#1}}

\newcommand{\annotation}[1]{\marginpar{\tiny #1}}

\newcommand\MD[1]{\annotation{\md{#1}}}

\begin{document}

\title{Lyapunov exponents and nonadapted measures for dispersing billiards}
\author{Vaughn Climenhaga\thanks{Department of Mathematics, University of Houston, Houston, TX 77204, USA. Email: climenha@math.uh.edu}
\and Mark F. Demers\thanks{Department of Mathematics, Fairfield University, Fairfield CT 06824, USA.  Email: mdemers@fairfield.edu} 
\and Yuri Lima\thanks{Departamento de Matem\'atica, Centro de Ci\^encias, Campus do Pici,
Universidade Federal do Cear\'a (UFC), Fortaleza -- CE, CEP 60455-760, Brasil. Email: yurilima@gmail.com}
\and Hongkun Zhang\thanks{Department of Mathematics and Statistics, University of Massachusetts Amherst, MA. Email: hongkunz@umass.edu}}

\date{\today}

\maketitle

\begin{abstract} 
For hyperbolic systems with singularities, such as dispersing billiards, Pesin theory as developed by Katok and Strelcyn applies to measures that are ``adapted'' in the sense that they do not give too much weight to neighborhoods of the singularity set. The zero-entropy measures supported on grazing periodic orbits are 
nonadapted, but it has been an open question whether there are nonadapted measures with positive entropy. We construct such measures for any dispersing billiard with a periodic orbit having a single grazing collision; we then use our construction to show that the thermodynamic formalism for such billiards has a phase transition even when one restricts attention to adapted or to positive entropy measures.
\end{abstract}

\section{Introduction}
\label{intro}

Sinai introduced dispersing billiards as a mathematically tractable example of a mechanical system
where the Boltzmann hypothesis can be verified \cite{Sinai}.  Since then,
mathematical billiards, and the Lorentz gas in particular, have become central models in 
mathematical physics.  Although the billiard map associated with a Lorentz gas has discontinuities,
it preserves a smooth invariant measure, which we denote by $\musrb$.  Classical results focus on 
establishing quantitative statistical properties such as the Central Limit Theorem, decay of correlations and large deviation estimates with respect to $\musrb$, and they use a variety of techniques, including
Markov partitions and sieves \cite{bsc1, bsc2}, Young towers \cite{Young}, spectral analysis of the
transfer operator \cite{dz} and, most recently, projective cones \cite{DL}.

In the past years, there has been significant interest in the existence and properties of other invariant measures
for the billiard map.  In particular, \cite{cwz, bd1, bd2} study the existence and uniqueness of 
equilibrium measures with respect to
the geometric family of potentials, $- t \log J^uT$, $t \ge 0$, where $J^uT$ denotes the unstable
Jacobian of the billiard map $T$.  A renewed interest in the construction of Markov partitions
for surface diffeomorphisms stemming from \cite{Sarig} has also led to symbolic codings 
for more general hyperbolic measures for billiards \cite{lima matheus,alp}.  
These constructions in turn lead to more questions regarding properties of these invariant measures,
such as possible Lyapunov spectra and phase transitions in the associated pressure functions.

For a class of planar dispersing billiard maps corresponding to a finite horizon Lorentz gas, we
present a construction that answers two open questions in the billiards literature.
The first relates to the construction of countable Markov partitions for
billiards, which relies on Pesin theory as developed by Katok and Strelcyn; this theory applies to invariant probability measures $\mu$ that are \emph{adapted}, meaning that
\begin{equation}\label{eqn:adapted}
\int |\log d(x,\mathcal{S})| \,d\mu < \infty,
\end{equation}
where $\mathcal{S}$ is the singular set for $T$.\footnote{See \cite[p.\ 342]{bH22} for the story of the genesis of this condition, or rather of the slightly stronger condition \eqref{a} that appears in \cite{KS86,yP92}. 
The term ``adapted'' 
was introduced in \cite{lima sarig}, in the context of smooth flows equipped with an invariant measure, as a condition on a Poincar\'e section (requiring that the boundary $\mathcal{S}$ satisfies \eqref{eqn:adapted} with respect to the induced measure).
It has subsequently been
used as a condition on measures in systems with singularities \cite{lima matheus, bd1}. 
For planar dispersing billiards, it is equivalent to assuming the finiteness of the positive Lyapunov exponent almost everywhere (see the proof of Theorem~\ref{thm-pressure}). }
Using this theory,
Lima and Matheus showed that if $\mu$ is adapted,
ergodic and has positive metric entropy, then there is a countable Markov partition for the billiard map $T$ which provides
a symbolic coding for $\mu$-almost every point \cite{lima matheus}. Therefore, a natural question 
arises: are there $T$-invariant {\em nonadapted} measures with positive entropy? We answer this question 
in the affirmative.

\begin{theorem}\label{thm-non-adapted}
There are finite horizon planar dispersing billiards with invariant probability measures
of positive entropy that are not adapted.
\end{theorem}

The second question is related to the behavior of the pressure function $P(t)$ for the family of geometric potentials
$- t \log J^uT$ for $t$ near 0, where $J^uT$ denotes the unstable Jacobian of $T$.   The pressure function is defined as
\begin{equation}
\label{eq:P}
P(t) = \sup_{\mu \in \cM} \left\{ h_\mu(T) - t \int \log J^uT \, d\mu \right\}  ,
\end{equation}
where $\cM$ is the set of $T$--invariant probability measures.
The function $P(t)$ is finite and decreasing for $t \ge 0$, and \cite{bd2} proved it is
analytic on an interval $(0, t_*)$, for some $t_*>1$ depending on the billiard table.
While $P(0) < \infty$ is the topological entropy of the map \cite{bd1}, it is clear that for
a table with a grazing periodic orbit, $P(t) = \infty$ for $t<0$ due to the invariant
measure supported on such an orbit. We arrive at the following natural question: does one still
see a jump in $P(t)$ at $t=0$ if one restricts the supremum in \eqref{eq:P} to a smaller
class of measures? For example, letting $N_\ve(\cS)$ denote the $\ve$--neighborhood of 
$\cS$, we can restrict the supremum to measures satisfying:
\begin{enumerate}[(a)]
\item\label{a} There are constants $C, \alpha>0$ such that $\mu(N_\ve(\cS)) \le C\ve^\alpha$ for all $\ve >0$ (such a measure must be adapted); or
\item $\mu$ has positive entropy.
\end{enumerate}
Our construction answers this question in the affirmative for both restricted suprema under 
condition (a) or (b) whenever a certain
type of periodic orbit exists on the billiard table.
This implies that, for some classes of dispersing billiards, there is a phase
transition at $t=0$.

\begin{theorem}\label{thm-pressure}
There are finite horizon planar dispersing billiards such that each of them possesses a sequence of periodic orbits $(p_n)$
with positive and finite Lyapunov exponents, but tending to infinity as $n \to \infty$.
As a consequence, for all $t<0$ we have
\begin{equation}
\label{eq:var a}
\sup_{\mu \in \cM} \left\{ h_\mu(T) - t \int \log J^uT \, d\mu : \exists C, \alpha > 0 \mbox{ such that }N_\ve(\cS) \le C\ve^\alpha,\forall \ve > 0 \right\} = \infty .
\end{equation}
Moreover, it follows from Theorem~\ref{thm-non-adapted} that for such dispersing billiards and for all $t < 0$
we have
\begin{equation}
\label{eq:var b}
\sup_{\mu \in \cM} \left\{ h_\mu(T) - t \int \log J^uT \, d\mu : \mbox{ $\mu$ has positive entropy } \right\} = \infty.
\end{equation}
\end{theorem}

The corresponding questions for the infinite horizon Lorentz gas have already been answered in the affirmative
in \cite{chernov troubetzkoy}.  There, the authors show that $P(0) = \infty$ and construct invariant measures with 
both infinite metric entropy as well as invariant measures with finite entropy, but infinite Lyapunov exponent.

\begin{remark}
We remark that for $t>0$, restricting to measures that satisfy both {\em (a)} and {\em (b)} does not change the pressure $P(t)$;
specifically, the equilibrium states produced in {\em \cite{bd2}} satisfy both conditions for all $t \in (0, t_*)$.
However, it is not known whether the measure of maximal entropy $\mu_0$ (in the case $t=0$) satisfies {\em (a)}.  
Indeed, {\em \cite{bd1}} produces the weaker bound $\mu_0(N_\ve(\cS)) \le C |\log \ve|^{-\gamma}$ for some $\gamma>1$; this is sufficient
for $\mu_0$ to be adapted and hyperbolic.

We state Theorem~\ref{thm-pressure} using the stronger condition {\em (a)} since it makes for a stronger statement:
even restricting to this smaller class of measures, the pressure is infinite when $t<0$.
\end{remark}

\subsection{Acknowledgements}

This work started as one of the projects proposed in the workshop
{\em Equilibrium states for dynamical systems arising from geometry}, held at the American Institute of Mathematics 
in July 2019. The authors are grateful to AIM for its hospitality.
VC was partially supported by NSF grants DMS-1554794 and DMS-2154378, and by a Simons Foundation Fellowship.
MD is partially supported by NSF grant DMS-2055070.
YL was supported by CNPq and Instituto Serrapilheira, grant ``Jangada Din\^{a}mica: Impulsionando Sistemas Din\^{a}micos na 
Regi\~{a}o Nordeste''. HZ is partially supported by NSF grant DMS-2220211 and Simons Foundation 706383.


\section{Setting}

Our billiard table is defined by placing finitely many pairwise 
disjoint closed, convex sets $O_i$, $i = 1, \ldots, d$, in $\mathbb{T}^2 = \mathbb{R}^2 / \mathbb{Z}^2$.  We assume
that the boundaries $\partial O_i$ are $C^3$ curves with strictly positive curvature.  The billiard table is then defined as
$\cQ := \mathbb{T}^2 \setminus (\bigcup_{i=1}^d O_i)$.  The billiard flow is induced by the motion of a point particle
traveling at unit speed in $\cQ$ and reflecting elastically at collisions with $\partial \cQ$.  The associated
billiard map, which we denote by $T$, is the Poincar\'e map under the flow with respect to $\partial \cQ$.

Parametrizing $\partial O_i$ for each $i$ by arclength $r$ (oriented clockwise), and letting $\vf \in [-\pi/2, \pi/2]$ denote the angle made
by the post-collision velocity vector and the outward normal to $\partial O_i$ at the point of collision, 
we represent the phase space of our billiard map
as $M = \bigcup_{i=1}^d (\partial O_i \times [-\pi/2, \pi/2])$.  We denote by $\tau(x)$ the flight time from $x = (r,\vf) \in M$ 
to $T(x)$.  We assume that the table has {\em finite horizon}, i.e. $\sup_{x\in M}\tau(x)<\infty$, is uniformly bounded. 
The finite horizon assumption together with the fact that the scatterers are disjoint implies that,
there exist constants $\tau_{\min}$, $\tau_{\max} \in \mathbb{R}$, such that 
$0 < \tau_{\min} \le \tau(x) \le \tau_{\max} < \infty$ for all $x \in M$.

It is well-known that $T$ preserves a smooth invariant probability measure $\musrb$, such that
$$
d\musrb(x) = \frac{1}{2|\partial Q|} \cos \vf \, dr d\vf.
$$
Yet $T$ has discontinuities created by grazing collisions.  Set 
$\cS_0 = \left\{ (r,\vf) \in M : \vf = \pm \frac{\pi}{2} \right\}$ and denote by
$\cS_n = \bigcup_{i=0}^n T^{-i}\cS_0$ and $\cS_{-n} = \bigcup_{i=0}^n T^i \cS_0$ the singularity sets 
for $T^n$ and $T^{-n}$, respectively.  Then $T^n$ is a $C^2$ diffeomorphism of $M \setminus \cS_n$ 
onto $M \setminus \cS_{-n}$.
Note that in this notation, $\cS_1$ is precisely the same as the set $\cS$ discussed in Section~\ref{intro}.


\subsection{Strategy of proof}

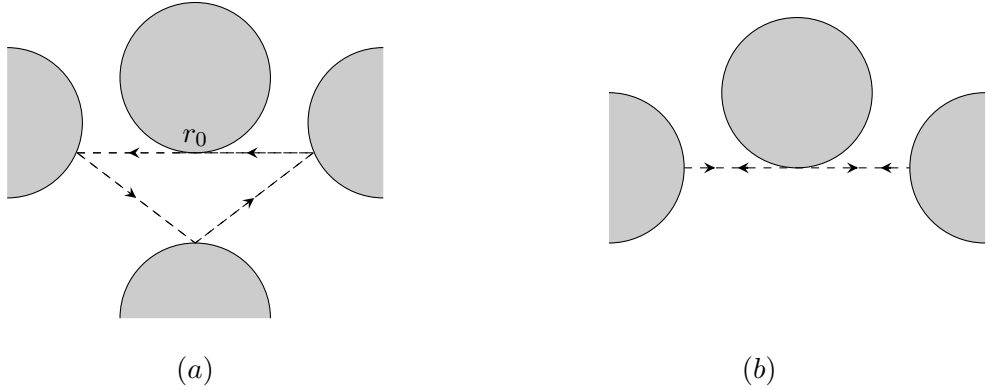
\begin{figure}[ht]
\begin{centering}
\begin{tikzpicture}[x=10mm,y=10mm]

 
 \filldraw[fill=black!20!white, draw=black] (4,1.6) arc (270:450:1);
\filldraw[fill=black!20!white, draw=black] (9,3.6) arc (90:270:1);
\filldraw[fill=black!20!white, draw=black] (7.5,0) arc (0:180:1);
\filldraw[fill=black!20!white, draw=black] (6.5,3.2) circle (1);

\draw[dashed] (4.93,2.2) -- (8.07,2.2) -- (6.5,1) -- (4.93,2.2);

\begin{scope}[dashed,decoration={
    markings,
    mark=at position 0.5 with {\arrow[thick]{stealth}}}
    ] 
    \draw[postaction={decorate}] (8.07,2.2)--(6.3,2.2);
    \draw[postaction={decorate}] (6.3,2.2)--(4.93,2.2);
    \draw[postaction={decorate}] (4.93,2.2)--(6.5,1);
    \draw[postaction={decorate}] (6.5,1)--(8.07,2.2);
\end{scope}

\node at (6.5,2.4){\small $r_0$};

\node at (6.5,-.7){\small$(a)$};


\filldraw[fill=black!20!white, draw=black] (12,1) arc (270:450:1);
\filldraw[fill=black!20!white, draw=black] (17,3) arc (90:270:1);
\filldraw[fill=black!20!white, draw=black] (14.5,3) circle (1);

\draw[dashed] (13,2) -- (16,2);
\draw[-stealth,thick] (13.4,2)--(13.41,2);
\draw[-stealth reversed, thick] (13.85,2)--(13.86,2);
\draw[-stealth, thick] (15.3,2) -- (15.31,2);
\draw[-stealth reversed, thick] (15.75,2) -- (15.76,2);

\node at (14,-.7){\small$(b)$};

 \end{tikzpicture}
\caption{(a) A periodic orbit of period 4 with grazing collision at $r_0$. The corresponding
periodic point is $x_0=(r_0,\pi/2)$. Since there is only one grazing collision, this configuration satisfies our assumption. 
(b) A periodic orbit of period 4 with two grazing collisions.
It does not satisfy our assumption. 
Both figures are local and are not meant to illustrate the finite horizon condition.}
\label{fig:orbits}
\end{centering}
\end{figure}

The class of billiards tables $\cQ$ we will consider admits the following type of periodic orbit with a grazing collision.
If we choose a direction and follow the orbit around one full cycle, then either all grazing collisions along the orbit 
occur on our right, or all grazing collisions along the orbit occur on our left.  See Figure~\ref{fig:orbits}(a) for 
an example of an orbit that satisfies our assumption, and Figure~\ref{fig:orbits}(b) for an example of an orbit that 
does not. As we will see in the next section, this assumption allows us to easily understand how a neighborhood 
in $\cS_0$ of the periodic point is cut under iterations.

The proofs of Theorems \ref{thm-non-adapted} and \ref{thm-pressure} will both 
use the existence of a Cantor rectangle defined in a one-sided neighborhood of the periodic point $x_0 \in \cS_0$,
and it will be divided into the following steps:
\begin{enumerate}[1.]
  \item  The point $x_0$ has a local stable manifold $W^s(x_0)$ and a local unstable manifold $W^u(x_0)$
  of positive length that terminate at $x_0$. Both manifolds are $C^{1+\rm{Lip}}$, and the wedge between $W^s(x_0)$ and $W^u(x_0)$ belongs to
  a single component of
  $M\setminus(\cS_p \cup \cS_{-p})$, where $p$ is the period of the orbit.
  \item Construct a Cantor rectangle $\cR$ formed by intersections of local stable and local unstable
  manifolds, with $W^s(x_0)$ and $W^u(x_0)$ comprising two boundaries of $\cR$, and such that
  $\musrb(\cR) > 0$.
  \item Use the product structure of the Cantor set and the mixing property of $\musrb$ to define a 
  horseshoe $K$ on $\cR$.
  \item Use the coding provided by the horseshoe $K$ to construct periodic orbits with large period and arbitrarily large Lyapunov exponent.
  \item Use the coding for the induced map on the horseshoe $K$ to define a Bernoulli measure whose projection to $M$ is nonadapted. 
\end{enumerate}
One of the difficulties involved in the construction is that $DT(x_0)$ is not defined, and $DT(x)$ converges
to infinity as $x$ converges to $x_0$. This means that $T$ has huge expansion near $x_0$, and it prevents
the use of classical techniques, such as the graph transform method. Nevertheless, the assumption on the periodic
orbit allows us to bypass this issue.

\medskip
\noindent
{\bf Notation.}  Throughout the text, we will use the notation $A \approx B$ if there exists
a constant $C>0$, depending only on the billiard table and the definition of the homogeneity strips,
such that $C^{-1} A \le B \le C A$.


\section{Construction of the Cantor rectangle $\cR$}

For ease of the presentation, we assume that the configuration of $x_0=(r_0,\pi/2)$ is the one depicted
in Figure~\ref{fig:orbits}(a), i.e. $T^4(x_0)=x_0$ with a single grazing collision in one full cycle.
Let $\mathcal K$ denote the curvature function
of the obstacles. We recall from \cite[Section~4.5]{chernov book} that $T$ has invariant cones:
\begin{align}\label{cones}
\begin{array}{l}
\mathcal C^s_x=\{(dr,d\vf)\in T_xM: -\mathcal K-\cos\vf/\tau\leq d\vf/dr\leq -\mathcal K\}\\
\\
\mathcal C^u_x=\{(dr,d\vf)\in T_xM: \mathcal K\leq d\vf/dr\leq \mathcal K+\cos\vf/\tau\}.
\end{array}
\end{align}
Since the table has finite horizon and the obstacles are disjoint and $C^3$ with strictly positive curvature, 
these cones are uniformly transverse and bounded away from the lines $r={\rm const}$ and $\vf={\rm const}$.
According to \cite[Eq. (4.19)]{chernov book}, vectors in these cones contract and expand uniformly:
there exists $C_e>0$ such that
\begin{equation}
\label{eq:hyp}
\| DT^n(x) \vec{v} \| \ge C_e \Lambda^n \| \vec{v}\| \, \text{ for all } \vec{v} \in \cC^u_x \text{ and all }n\geq 1 ,
\end{equation}
where $\Lambda = 1 + 2 \cK_{\min} \tau_{\min}$.  Similar bounds hold for $\vec{v} \in \cC^s_x$.

We call a $C^1$ curve $W$ a stable curve if the tangent vector at each $x \in W$ lies in $\cC^s_x$.
Unstable curves are defined similarly.


\subsection{Existence of stable and unstable manifolds for $x_0$}

In this section, we will construct one-sided stable and unstable manifolds for $x_0$.  The main result
of this section is the following proposition.

\begin{proposition}
\label{prop:unstable}
Let $x_0 \in \cS_0$ be a periodic point with grazing collision as described above.  Then $x_0$ has one-sided local stable and 
unstable manifolds, which are $C^{1+{\rm Lip}}$ curves in $M$. 
\end{proposition}

Proposition~\ref{prop:unstable} follows from Lemma~\ref{lemma-W^u}.

\begin{figure}[ht]
\begin{centering}
\begin{tikzpicture}[x=7mm,y=7mm]

\draw[thick] (0,0) -- (14,0);
\draw[thick] (7,0) to[out=235, in=45] (4.9,-2.3) to[out=225, in=35] (3, -4);
\draw (7,0) to[out=240, in=60] (6.5,-1) to[out=240, in=60] (4.5,-4);

\draw[thick] (7,0) to[out=-45, in=135] (9.1,-2.3) to[out=-45, in=145] (11, -4);
\draw (7,0) to[out=-50, in=120] (7.5,-1) to[out=300, in=120] (9.5,-4);

\node at (-.4,0){\small $\cS_0$};
\node at (4.9,.4){\small $I_-$};
\node at (8.9,.4){\small $I_+$};
\node at (6.9,.4){\small $x_0$};
\node at (2,-4){\small $T^4I_-$};
\node at (4.3,-4.5){\small $W^u(x_0)$};
\node at (9.3, -4.5){\small $W^s(x_0)$};
\node at (12.2, -4){\small $T^{-4}I_+$};
\node at (6.9, -3){\small $\Delta$};

 \end{tikzpicture}
\caption{A $\delta_0$-neighborhood of $x_0$.}
\label{fig:neigh}
\end{centering}
\end{figure}
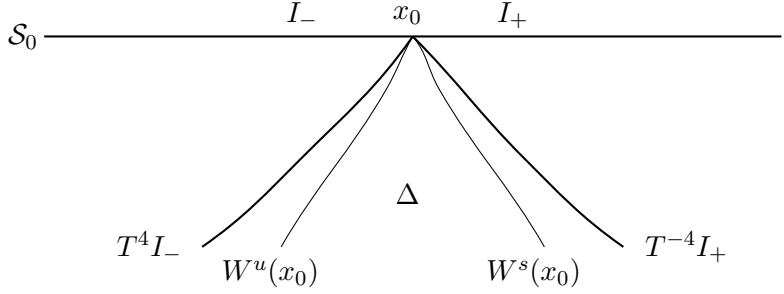

Let $O_1$ be the obstacle of the grazing collision and let $M_1$ denote the corresponding 
component of the phase space.
We assume that the orientation of the trajectory is counter-clockwise, and that the parametrization
of the obstacles is clockwise.
Consider a neighborhood $I=[r_1,r_2]\times\{\pi/2\}\subset \cS_0$ of $x_0$ and let 
$I_-=[r_1,r_0]\times\{\pi/2\}$ and $I_+=(r_0,r_2]\times\{\pi/2\}$.
Since $\cS_0$ is a dispersing wave front, 
we may choose $r_1$ and $r_2$ sufficiently close to $r_0$ such that the following behavior holds for the first four iterates of $I$:
\begin{enumerate}[$\circ$]
\item If $x\in I_- \setminus \{ x_0 \}$ then $r(T(x))>r(T(x_0))$, $r(T^2(x))<r(T^2(x_0))$, $r(T^3(x))>r(T^3(x_0))$,
so that $T^4(x)$ hits $O_1$ with  $r(T^4(x))<r(T^4(x_0))=r_0$. 
\item If $x\in I_+$ then $r(T(x))<r(T(x_0))$, $r(T^2(x))>r(T^2(x_0))$, $r(T^3(x))<r(T^3(x_0))$,
so that $T^4(x)$ {\em does not} hit $O_1$ in a neighborhood $N_\ve(x_0)$ of $x_0$. 
\end{enumerate}
This means that $I$ gets cut after four iterations: $T^4(I_-)$ is a curve in $M_1$ terminating at $x_0$,
and $T^4(I_+)\cap N_\ve(x_0) = \emptyset$ for some $\ve >0$. 
Analogous facts hold for $I_+$ under $T^{-1}$.  
Thus we may choose $r_1$ and $r_2$
sufficiently close to $r_0$ so that $T^4$ is continuous on $I_-$ and $T^{-4}$ is continuous on $I_+$. 
See Figure~\ref{fig:neigh}.

Consider the sector $\cA$ defined by rotating the horizontal segment
$I_- = \{ (r, \vf) : \vf = \pi/2, r \in [r_1, r_0] \}$ around $x_0$ by an angle of $\pi/2$ radians,
until $I_-$ reaches the vertical segment $V = \{ (r, \vf) : r = r_0, \vf \in [\pi/2 - (r_0-r_1), \pi/2] \}$.  
Both segments $I_-$ and $V$ terminate on $x_0$, a fixed point for $T^4$.   
Moreover, both $T^4(I_-)$ and $T^4(V)$ are unstable curves, tangent at $x_0$
(the cone at $x_0$ consists of a single line with slope $\mathcal K(r_0)$, see (\ref{cones}),
and contained strictly between
$I_-$ and $V$.  Since $T$ is orientation preserving, $T^4(I_-)$ comprises the upper boundary of
$T^4(\cA)$ and $T^4(V)$ comprises its lower boundary.
Repeating this argument, since the sector $\cA$ is strictly contracted by $T^4$,
the intersection $T^8(\cA) \cap S$, which has upper boundary $T^8(I_-) \cap \cA$ and lower boundary 
$T^8(V) \cap \cA$, lies strictly inside the region with upper boundary $T^4(I_-)$ and lower boundary $T^4(V)$.

Inductively, we see that for each $n \ge 1$ the curve $T^{4n}(I_-)$ contains an unstable curve 
$\gamma_n(x_0)$, defined by an equation
$\{(r,\vf_n(r)):r\in [r_1,r_0]\}$ where $\vf_n:[r_1,r_0]\to\R$ is an increasing $C^2$ function 
with $\vf_n(r_0)=\pi/2$ and $\vf_n'(r_0)=\mathcal K(r_0)$.
Moreover, the sequence $(\vf_n)$ is decreasing, and $(\vf_n')$ is bounded since it must belong
to the unstable cone for each $n$. 
Thus, $(\vf_n)$ converges to a function $\vf^u:[r_1,r_0]\to \R$. This function defines an unstable curve, 
which indeed is the local unstable manifold at $x_0$, as we will prove below.
We need the following auxiliary lemma.

\begin{lemma}\label{lemma-growth}
Suppose $n \ge 1$ and $W$ is an unstable curve passing through $x_0$
such that $T^{-4k}W$ is an unstable curve for all $1 \le k \le n$. The following holds for all $z \in W$:
\begin{enumerate}[i.]
  \item[{\rm (a)}] $d(T^{-4}z ,x_0)\approx d(z,x_0)^2$;
  \item[{\rm (b)}] $d(T^{-4n} z, x_0) \le C_0 d(z, x_0)^{(3/2)^n}$, for some $C_0>0$ independent of $W$ and $n$.
\end{enumerate}
The analogous statement holds for stable curves $W$ passing through $x_0$ such that
  $T^{4k}W$ is a stable curve for all $1 \le k \le n$.
\end{lemma}


This means that near $\cS_0$, the map $T^4$ exhibits a huge expansion along unstable
curves.

\begin{proof}
We consider the homogeneity strips 
$$
H_k=\left\{x=(r,\varphi)\in M: \varphi\in \left[\tfrac{\pi}{2}-k^{-2},\tfrac{\pi}{2}-(k+1)^{-2}\right)\right\},\ k\geq 1.
$$
Notice that $\cos\vf(x)\approx k^{-2}$ for $x\in H_k$.
Let $z$ and $W$ be as in the statement. Let $V $ denote the subcurve of $W$
with endpoints $z$ and $x_0$, and write $V=\bigcup_{k\geq j}V_k$, where $V_k=V \cap H_k$. 
Since $V_k$ is uniformly transversal to the
boundaries of $H_k$, we have $|V_k|\approx k^{-3}$, hence
$$
|V|=\sum_{k\geq j} |V_k|\approx \sum_{k\geq j}k^{-3}\approx j^{-2}.
$$
We will show that $|T^{-4}V|\approx j^{-4}$, which will complete the proof of part (a) the lemma.

By \cite[Equation (4.20)]{chernov book}, for all $x \in H_k$ we have
$$
J^uT^4(T^{-4}x) \approx \frac{1}{\cos\vf(x)}   \approx k^{2}.
$$
Hence $|T^{-4}V_k|\approx |V_k|\cdot k^{-2}\approx k^{-5}$ and so
$$
|T^{-4}V|=\sum_{k\geq j}|T^{-4}V_k|\approx \sum_{k\geq j}k^{-5}\approx j^{-4},
$$
which concludes the proof of part (a).

To prove (b), note that (a) yields a constant $C > 0$ such that $d(T^{-4}z, x_0) \le C d(z, x_0)^2$.  Without loss of generality, we may assume that $d(z, x_0) \le C^{-2}$, as this must 
occur in a fixed number of iterates due to the uniform expansion given by \eqref{eq:hyp}.
Statement (b) then follows directly by iterating (a) since for such $z$ we have $C d( T^{-kn}z, x_0)^2
\le d(T^{-kn}z, x_0)^{3/2}$.
\end{proof}

%
%
%

\begin{lemma}\label{lemma-W^u}
The curve
$$
W^u(x_0):=\{(r,\vf^u(r)):r\in [r_1,r_0]\}
$$
is the local unstable manifold at $x_0$ for $T^4$, i.e.
$d(T^{-4n}x,x_0)\xrightarrow[]{n\to\infty} 0$ for every $x\in W^u(x_0)$.
In particular, $W^u(x_0)$ is a $C^{1+{\rm Lip}}$ curve.
\end{lemma}

By time reversibility of billiard systems, $x_0$ also has a $C^{1+{\rm Lip}}$ local stable manifold $W^s(x_0)$, see Figure~\ref{fig:neigh}.

\begin{proof}
It is clear that the function $\vf^u$ is Lipschitz and that $W^u(x_0)$ is an unstable curve. 
Fix $x\in W^u(x_0)$, and consider a stable curve $W$ passing through $x$.
Fix $n\geq 1$. For $1 \le k \le n$, recall that $\gamma_k \subset T^{4k}I_-$
denotes the graph of $\vf_k$ over the interval $[r_1, r_0]$.
Let $y$ be the intersection point of $W$ and $\gamma_n$. 
Consider the pre-images $T^{-4n}x$ and $T^{-4n}y$, and consider the triangle $\alpha_1\alpha_2\alpha_3$
formed by the three curves:
\begin{align*}
\alpha_1&= \text{segment of $\cS_0$ joining $x_0$ and $T^{-4n}y$}\\
\alpha_2&= \text{curve of $T^{-4n}W$ joining $T^{-4n}y$ and $T^{-4n}x$}\\
\alpha_3&= \text{curve of $W^u(x_0)$ joining $T^{-4n}x$ and $x_0$.}
\end{align*}
Since $\alpha_2$ is a stable curve and $\alpha_3$ is an unstable curve, and cones
make angles bounded away from zero with respect to horizontal and vertical curves of $M$, it follows that 
$d(T^{-4n}x,x_0)\leq {\rm const}\times d(T^{-4n}y,x_0)$.
Since $\gamma_k$ is an unstable curve for each $k \ge 1$ and $T^{4(k-n)}y \in \gamma_k$,
by Lemma \ref{lemma-growth}(b) we have
$$ 
d(T^{-4n}y,x_0) \le C d(T^{-4(n-1)}y,x_0) \le CC_0 d(y,x_0)^{(3/2)^{n-1}} , 
$$
and so $d(T^{-4n}x,x_0)\leq {\rm const}\times d(y,x_0)^{(3/2)^{n-1}}$,
which converges to $0$
as $n\to\infty$.

It remains to prove that $\vf^u$ is $C^{1+{\rm Lip}}$.  According to the first part of the proof, for each 
$x\in W^u(x_0)\backslash \{x_0\}$, the curve $W^u(x_0)$ is the unstable
manifold of $x$. By \cite[Corollary~4.61]{chernov book},
$W^u(x_0)$ is $C^1$ and its derivative is Lipschitz.
\end{proof}

\subsection{The Cantor rectangle $\cR$}

In the last section, we showed the existence of $C^{1+{\rm Lip}}$ local stable/unstable manifolds
$W^{s/u}(x_0)$ at $x_0$
for the map $T^4$. Note that, since these are a stable/unstable curve respectively, they make an angle 
at $x_0$ that is in\footnote{Indeed, since the slopes are $-\cK(x_0)$ and $\cK(x_0)$, respectively,
they form an angle equal to 2Arctan$(1/\cK(x_0))$.} 
$(0,\pi/2)$. We denote by $\Delta$ the region between these two curves, see Figure~\ref{fig:neigh}.
In this section, we will construct a Cantor rectangle of stable and unstable manifolds with $x_0$ as one corner point
of the rectangle.

First we recall some terminology.  We call $D$ a solid rectangle in $M$ if $\partial D$ comprises
four smooth curves,
two local stable and two local unstable manifolds. 
Given a solid rectangle $D$, let $\mathfrak{S}^s(D)$ and
$\mathfrak{S}^u(D)$ denote the set of local stable and unstable manifolds, respectively, of points in $D$ that
do not terminate in the interior of $D$.
Define the Cantor rectangle 
\[
R(D) = \mathfrak{S}^s(D) \cap \mathfrak{S}^u(D) \cap D.  
\]
By construction, $R(D)$ has a locally maximal hyperbolic product structure.  
We will work with such
Cantor rectangles $R$, closed sets formed by locally maximal intersections of stable and unstable manifolds such that
Leb$(R)>0$.  Conversely, given such a Cantor rectangle
$R$, we denote by $D(R)$ the smallest solid rectangle containing $R$.

Unlike \cite[Section~7.11]{chernov book}, we do not restrict ourselves to $H$--manifolds, i.e. to local stable (respectively unstable) manifolds whose forward (respectively backward) trajectories lie in a single homogeneity strip at each step.
We do this because the rectangle we construct, $\cR$, will necessarily cross infinitely many homogeneity strips
in order to incorporate $x_0$ in its boundary.

Now consider a segment of $W^u(x_0)$ that lies between homogeneity strips of index $k_1$ and $k_1^2$.  Call this curve
$A_1$. Since $A_1$ is a finite union of homogeneous unstable manifolds, almost every point on $A_1$ has stable
manifolds of positive length \cite[Theorem~5.70]{chernov book}.  

Now consider $T^{-4}A_1 \subset W^u(x_0)$. If $W$ is a local stable manifold of a point in $A_1$,
then $T^{-4}W$ is a local stable manifold of a point in $T^{-4}A_1$, and it has one connected component 
intersecting $\Delta$.
This component cannot be cut until it exits a $\delta_0$--neighborhood of $x_0$,
and it cannot cross $W^s(x_0)$ (since stable
manifolds do not intersect).  Therefore,
\[
|W^s(y)| \ge \min \left\{ \frac{|W^s(T^4y)|}{ J^sT^4(y)}, \delta_0 \right\} \, \mbox{ for all $y \in T^{-4}A_1$}.
\] 
Note that $J^sT^4(y) \approx \cos \vf(y)$ so that the expansion is large in a $\delta_0$--neighborhood of $x_0$.  

Applying this construction inductively, we are guaranteed the existence of increasingly longer stable manifolds
for points in $W^u(x_0)$ as we approach $x_0$, and so a positive measure of such manifolds 
extend a minimum distance of $\delta_1 > 0$ inside $\Delta$.

By time-reversal symmetry, this construction also implies that almost every point on $W^s(x_0)$ has an unstable
manifold of positive length, and that unstable manifolds of a minimum length $\delta_1$ comprise a positive measure set
of points inside $W^s(x_0)$.

Fixing once and for all $\delta_0,\delta_1 >0$ small enough such that the above constructions hold in the 
$\delta_0$--neighborhood of $x_0$, we obtain a (locally maximal) closed, Cantor rectangle $\cR$ 
comprised of the intersection of a positive measure set of stable and unstable manifolds, having a segment
of $W^s(x_0)$ along one boundary, a segment of $W^u(x_0)$ along another, and containing the 
point $x_0$ as one of its corner points.


\subsection{Definition of the horseshoe}

We let $G=T^4$, and construct a horseshoe for a power of $G$.
The idea is to have two branches for the horseshoe, which are two stable rectangles in $\cR$, one of them
containing $x_0$. We first show that an iterate of $W^u(x_0)$ crosses $W^s(x_0)$, i.e. there is an 
homoclinic intersection. This is a consequence of standard results in the theory of dispersing billiards.

We need a little more terminology.  Given a Cantor rectangle $R$, we call $S \subset R$ an $s$--subrectangle
of $R$ if for each $x \in S$, $W^s(x) \cap S = W^s(x) \cap R$.  Similarly, $U \subset R$ is a $u$--subrectangle 
if $W^u(x) \cap U = W^u(x) \cap R$ for all $x \in U$.  We say a local unstable manifold $W^u$
{\em fully crosses} $R$ if $W^u \cap \mathring{D}(R) \ne \emptyset$ and $W^u$ does not terminate in the
interior of $D(R)$.  A Cantor rectangle $R'$ $u$--crosses $R$ if every unstable manifold $W^u \in \mathfrak{S}^u(R')$
fully crosses $R$.  Analogous definitions hold for stable manifolds and $s$--crossings.

Given a Cantor rectangle $R$, $T^n(R)$ is a finite union of (maximal) Cantor rectangles (recall that stable manifolds cannot
be cut by singularities of $T$).  We label them $R_{n,i}$.  Then each $T^{-n}(R_{n,i})$ is an $s$--subrectangle
of $R$. 

With these preliminaries, we are ready to proceed with the construction of our horseshoe.

\begin{enumerate}[$\circ$]
\item By \cite[Lemma 7.87]{chernov book}, there is a finite collection of rectangles $R_1,\ldots,R_N$
of positive measure each such that any stable and unstable curve of length at least $\delta_1/2$ 
properly crosses at least one of the rectangles. Without loss of generality, we can assume that 
$W^u(x_0)$ properly crosses $R_1$ and $W^s(x_0)$ properly crosses $R_2$.
\item By  \cite[Lemma 7.90]{chernov book}, there is a `magnet' rectangle $R^*$ of positive measure, 
and a `high density' subset $\mathfrak{P}^* \subset R^*$, satisfying
the following property: 
if $R_{k,n,i}$ is a maximal rectangle in $T^n(R_k)$ and
$R_{k,n,i} \cap \mathfrak{P}^* \neq\emptyset$ where $n$ is large enough,
then $R_{k,n,i}$ $u$--crosses $R^*$.\footnote{Indeed, a slightly stronger property holds:
 for every $x\in R_{k,n,i}$ its local unstable manifold
$W^u(x)$ properly crosses $R^*$, a proper crossing being a full crossing whose distance from the
unstable boundary of $R^*$ is at least a fixed fraction of the unstable diameter of $R^*$.}
The analogous properties hold for maximal
rectangles $R_{k,-n,i}$ of $T^{-n}(R_k)$ and $s$--crossings of $R^*$.

\item Since $\musrb$ is mixing, there are $m,n>0$ such that $G^n(R_1)\cap \mathfrak{P}^*\neq\emptyset$
and $G^{-m}(R_2)\cap \mathfrak{P}^*\neq\emptyset$. Therefore there are an $s$-subrectangle $R'\subset R_1$
and a $u$-subrectangle $R''\subset R_2$ such that $G^n(R')$ $u$-crosses $R^*$ and 
$G^{-m}(R'')$ $s$-crosses $R^*$. 
This implies that $G^n(W^u(x_0))\cap G^{-m}(W^s(x_0))\neq\emptyset$, proving the desired
transverse intersection.
\end{enumerate}

Let $x_1\in W^u(x_0)\cap G^{-(m+n)}W^s(x_0)$. Take points $z\in W^u(x_0)\cap \cR$ and
$w\in W^s(x_0)\cap \cR$ such that $x_1$ belongs to the segment of $W^u(x_0)$ joining $x_0$ and $z$.
The points $x_0,z,[z,w],w$ define a solid rectangle $D$. Let $U_0,U_1$ be
disjoint  $s$-subrectangles
of $D$ such that $U_0$ contains $x_0$ as a vertex and 
the solid rectangle of $U_1$ contains $x_1$ in the interior of its $u$--side. 
If $\ell>0$ is large, then $G^\ell(U_0)$ $u$-crosses $D$. Also,
by construction, $G^{m+n}(U_1)$ intersects $W^s(x_0)$ transversally, and
so if $\ell_0>0$ is large enough
then $G^{\ell_0+m+n}(U_0)$ also $u$-crosses $D$. Hence, we can fix $\ell_0>0$ such that
for $\ell=\ell_0+m+n$ the sets $G^\ell(U_0),G^\ell(U_1)$ both $u$-cross $D$. 
Letting $f=G^{\ell}$, we obtain a horseshoe for $f$ as the intersection
\[
K := \bigcap_{n=-\infty}^\infty f^n(U_0 \cup U_1) \, .
\] 
By construction, $K$ is conjugate to a full shift in two symbols, i.e.
if $(\Sigma,\sigma)$ is the topological Markov shift with $\Sigma=\{0,1\}^{\mathbb Z}$,
then there is a measurable bijection $\pi:\Sigma\to K$ such that $\pi\circ\sigma=f\circ\pi$.

%
%

\section{Proof of Theorem \ref{thm-non-adapted}}
The idea to prove Theorem \ref{thm-non-adapted} is the following: points that spend many iterates near $x_0$ approach $x_0$ super-exponentially fast. 
This is expressed for unstable curves (and analogously for stable curves) by Lemma \ref{lemma-growth}.
Given $x\in\Delta$, let $W^s,W^u$ be a stable, unstable curve passing through $x$ such that
 $\{x_s\}=W^u\cap W^s(x_0)$  and $\{x_u\}=W^s\cap W^u(x_0)$ are defined. Hence
$$
d(x,x_0)\approx \max\{ d(x_s,x_0),d(x_u,x_0)\}.
$$
 

Denote an element of $\Sigma$ by $\un v=(v_n)_{n\in\mathbb Z}$.
Let $C=\{\un v\in\Sigma:v_0=1\}$ and $C_n=\{\un v\in C:v_1=\cdots=v_{n-1}=0\text{ and }v_n=1\}$
for $n\geq 1$.
Clearly $C=\bigcup_{n\geq 1}C_n$ modulo $\musrb$, and the first return map $\widetilde \sigma:C\to C$ of $\sigma$
to $C$ is a full topological Markov shift in an infinite countable alphabet with
$\widetilde \sigma=\sigma^n$ on $C_n$. For $n \ge 1$, let $x=\pi(\un v)$ with $\un v\in C_n$.

As in the previous paragraph, let $W^u$ be an unstable curve connecting 
$f(x)$ and $W^s(x_0)$ and
define $\{ x^1_s \} = W^u \cap W^s(x_0)$.  Let $W^s$ be a stable curve connecting $f^{n-1}(x)$
and $W^u(x_0)$ and define $\{ x^{n-1}_u \} = W^s \cap W^u(x_0)$. 
Since $f^k(x^1_s) \in W^s(x_0)$ and $f^{-k}(x^{n-1}_u) \in W^u(x_0)$ for $k \ge 0$, we may
apply Lemma \ref{lemma-growth}(b) and obtain that:
\begin{enumerate}[$\circ$]
\item $d(f^{\lfloor \frac{n}{2}\rfloor}x^1_s,x_0) \le C_0 d(x^1_s,x_0)^{(3/2)^{\left\lfloor \frac{n}{2}\right\rfloor}}$.
\item $d(f^{-\lfloor \frac{n}{2}\rfloor}x^{n-1}_u,x_0) \le C_0 d(x^{n-1}_u,x_0)^{(3/2)^{\left\lfloor \frac{n}{2}\right\rfloor}}$.
\end{enumerate}
Since $W^s$ is a stable curve, $f^{-\lfloor \frac{n}{2} \rfloor}(W^s)$ is again a stable curve connecting
$f^{\lfloor \frac n2 \rfloor}(x)$ with $W^u(x_0)$.  Similarly, since $W^u$ is an unstable curve,
$f^{\lfloor \frac n2 \rfloor}(W^u)$ is an unstable curve connecting $f^{\lfloor \frac n2 \rfloor}(x)$ with
$W^s(x_0)$. 
 In particular, letting $\delta={\rm diam}(D)<1$,
there is a constant $C \ge 1$ such that 
\begin{equation}
\label{eq:dist}
d(f^{\lfloor \frac{n}{2}\rfloor}x, \cS_0)\leq C\delta^{(3/2)^{\left\lfloor \frac{n}{2}\right\rfloor}},\ \forall x\in \pi(C_n).
\end{equation} 

Let $b=\sum_{n\geq 1} \left(\frac 32\right)^{-{\lfloor \frac{n}{2}\rfloor}}<\infty$.
Take $p_n=b^{-1} \left(\frac 32\right)^{-{\lfloor \frac{n}{2}\rfloor}}$ for $n\geq 1$, and consider a Bernoulli probability measure $\widetilde\nu$ on
$C\cong \mathbb N^{\mathbb Z}$ such that $\widetilde \nu(C_n)=p_n$. 
This measure
induces a $\sigma$--invariant probability measure $\nu$ on $\Sigma$, which descends by $\pi$
to a $f$--invariant probability measure $\eta$ on $K$. 
Note that $\nu$ is supported on the
disjoint union $\bigcup_{n\geq 1}\bigcup_{0\leq k<n}\sigma^{k}C_n$, and 
$\eta$ is supported on the disjoint union $\bigcup_{n\geq 1}\bigcup_{0\leq k<n}f^{k}\pi(C_n)$.
Applying \eqref{eq:dist},
\begin{align*}
&\int |\log d(x, \cS_0)| \, d\eta\geq \int_{\bigcup_{n\geq 1}f^{\lfloor \frac{n}{2}\rfloor}\pi(C_n)} | \log d(x, \cS_0) | \, d\eta
= \int_{\bigcup_{n\geq 1}\pi(C_n)} | \log d(f^{\lfloor \frac{n}{2}\rfloor}x, \cS_0) | \, d\eta\\
&\geq \sum_{n\geq 1} \left| \log \left(C\delta^{(3/2)^{\left\lfloor \frac{n}{2}\right\rfloor}}\right) \right| \eta[\pi(C_n)] \ge
- \log C+
|\log \delta| \sum_{n\geq 1}p_n \left(\tfrac 32\right)^{\lfloor \frac{n}{2}\rfloor}=\infty.
\end{align*}

Now we estimate $h_\eta(f)$. 
Let $n$ be the return time function defining $\widetilde \sigma$. By the Abramov formula,
and using that $\widetilde\nu$ is a Bernoulli measure, we have
\begin{align*}
&\ h_\eta(f)=h_\nu(\sigma)=
\frac{h_{\widetilde\nu}(\widetilde \sigma)}{\int_C n(\underline v)d\widetilde \nu(\underline v)}
=\frac{\displaystyle\sum_{n=1}^\infty - p_n \log p_n}{\displaystyle\sum_{n=1}^\infty n p_n }
\approx
\frac{\displaystyle\sum_{n=1}^\infty \left(\tfrac 32\right)^{-\lfloor \frac{n}{2}\rfloor}{\left\lfloor \tfrac{n}{2}\right\rfloor} 
\log (3/2)}{\displaystyle\sum_{n=1}^\infty n \left(\tfrac 32\right)^{-\lfloor \frac{n}{2}\rfloor} }\\
&\approx\frac{\log (3/2)}{2}>0.
\end{align*}
Finally, let $\mu=\frac{1}{4\ell}\sum\limits_{k=0}^{4\ell-1}\eta\circ T^{-k}$, which is a $T$--invariant probability measure.
We claim that $\mu$ is a measure satisfying Theorem \ref{thm-non-adapted}. Firstly, we have
\begin{equation}
\label{eq:infinite}
\int  | \log d(x, \cS_0) | \, d\mu =\frac{1}{4\ell}\sum_{k=0}^{4\ell-1}\int | \log d(x, \cS_0) | \, d(\eta\circ T^{-k})\geq
\frac{1}{4\ell} \int | \log d(x, \cS_0) | \, d\eta=\infty.
\end{equation}
Since $\cS_0 \subset \cS_1$, this implies that the measure $\mu$ is not adapted.

Now, using that each $\eta\circ T^{-k}$ is $f$--invariant, it follows again by the Abramov formula
and the linearity of the metric entropy that
$$
h_\mu(T)=\frac{1}{4\ell}h_\mu(f)=\frac{1}{16\ell^2}\sum_{k=0}^{4\ell-1}h_{\eta\circ T^{-k}}(f)\geq \frac{\log (3/2)}{32\ell^2}>0.
$$
This concludes the proof of Theorem \ref{thm-non-adapted}.

\section{Proof of Theorem \ref{thm-pressure}}

The subset $\pi(C_n)$ contains a periodic orbit of period $n$ for $f$, equal to $y_n=\pi(\un v)$
where $v_0=1$ and $v_1=\cdots=v_{n-1}=0$.
Due to \eqref{eq:dist}, we have 
$d(f^{\lfloor \frac{n}{2}\rfloor}y_n, \cS_0)\leq C\delta^{(3/2)^{\left\lfloor \frac{n}{2}\right\rfloor}}$ and so,
using again \cite[eq.~4.20]{chernov book},
\[
J^uf(f^{\lfloor \frac{n}{2}\rfloor-1}y_n) \geq (\cos \vf(f^{\lfloor \frac n2 \rfloor} y_n))^{-1} \approx d(f^{\lfloor \frac{n}{2}\rfloor}y_n, \cS_0)^{-1} 
\ge C^{-1} \delta^{-(3/2)^{\left\lfloor \frac{n}{2}\right\rfloor}}.
\]
The Lyapunov exponent for the point $y_n$ with respect to $f$ is at least
\[
\begin{split}
\lim_{j \to \infty} \frac{1}{jn} \log J^uf^{jn}(y_n) & \ge \lim_{j \to \infty}  \frac{1}{jn} j \log J^uf(f^{\lfloor \frac{n}{2}\rfloor -1 }y_n) \\
& \geq \frac{1}{n} \log \left(C^{-1} \delta^{-(3/2)^{\left\lfloor \frac{n}{2}\right\rfloor}}\right)
\approx \frac{\left(\tfrac 32\right)^{\lfloor \frac{n}{2}\rfloor}|\log\delta|}{n} . 
\end{split}
\]
Therefore, The Lyapunov exponent for the point $y_n$ with respect to $T$ is at least
of the order of $ \left(\frac 32\right)^{\lfloor \frac{n}{2}\rfloor}/ 4\ell n$.
We have thus constructed a sequence of periodic orbits with finite, but arbitrarily large Lyapunov exponents.
This proves the first statement of Theorem~\ref{thm-pressure}.

Equation \eqref{eq:var a} follows immediately since each of the above periodic orbits contains no grazing collisions
and therefore each corresponding invariant measure is supported outside of an $\ve$--neighborhood 
of $\cS_0$ for some $\ve>0$.

To prove the final statement of the theorem, equation~\eqref{eq:var b}, let $\mu$ be the $T$--invariant
measure constructed in the proof of Theorem~\ref{thm-non-adapted}.   Using again \cite[(4.20)]{chernov book},
or \cite[eq.~(5.36)]{chernov book}, we have
$J^uT(x) \approx (\cos \vf (Tx))^{-1}$.  By the invariance of $\mu$ and estimate \eqref{eq:infinite},
we get that
\[
\int \log J^uT \, d\mu \approx \int - \log \cos (\vf \circ T) \, d\mu
= \int - \log \cos \vf \, d\mu \ge \int |\log d(x, \cS_0)| \, d\mu(x) = \infty \, .
\]
Since $h_\mu(T) >0$, this proves \eqref{eq:var b}.

\end{document}